\documentclass{amsbook}
\include {mak}
\begin{document}


\newcommand{\p}{\partial}
\newcommand{\eu}{{\rm e}}
\newcommand{\dd}{{\rm d}}

\newcommand{\bd}{\begin{definition}}                
\newcommand{\ed}{\end{definition}}                  
\newcommand{\bc}{\begin{corollary}}                 
\newcommand{\ec}{\end{corollary}}                   
\newcommand{\bl}{\begin{lemma}}                     
\newcommand{\el}{\end{lemma}}                       
\newcommand{\bp}{\begin{proposition}}            
\newcommand{\ep}{\end{proposition}}                
\newcommand{\bere}{\begin{remark}}                  
\newcommand{\ere}{\end{remark}}                     

\newcommand{\bt}{\begin{theorem}}
\newcommand{\et}{\end{theorem}}

\newcommand{\be}{\begin{equation}}
\newcommand{\ee}{\end{equation}}

\newcommand{\bit}{\begin{itemize}}
\newcommand{\eit}{\end{itemize}}
\newtheorem{theorem}{Theorem}[section]
\newtheorem{corollary}[theorem]{Corollary}
\newtheorem{lemma}[theorem]{Lemma}
\newtheorem{proposition}[theorem]{Proposition}
\theoremstyle{definition}
\newtheorem{definition}[theorem]{Definition}
\theoremstyle{remark}
\newtheorem{remark}[theorem]{Remark}
\newtheorem{example}[theorem]{Example}


\oddsidemargin 16.5mm \evensidemargin 16.5mm

\thispagestyle{plain}

\noindent {\small\sc Applicable Analysis and Discrete Mathematics}

\noindent {\scriptsize xx, No x (xxxx), x--x.}

\noindent{\scriptsize Available electronically at {\bf http://pefmath.etf.bg.ac.yu}}

\vspace{5cc}
\begin{center}

{\Large\bf  AN EQUIVALENT FORM OF YOUNG'S INEQUALITY WITH UPPER
BOUND
\rule{0mm}{6mm}\renewcommand{\thefootnote}{}
\footnotetext{
Keywords and Prases: Young's inequality, Legendre transform. }}

\vspace{1cc} {\large\it E. Minguzzi}

\vspace{1cc}
\parbox{24cc}{{\scriptsize\bf

%
Young's  integral inequality is complemented with an upper bound to
the remainder. The new inequality turns out to be equivalent to
Young's inequality, and the cases in which the equality holds become
particularly transparent in the new formulation.


}}
\end{center}

\vspace{1cc}

%
%


\section{Formulation of the theorem}

Let $\phi: [\alpha_1,\alpha_2]\rightarrow [\beta_1,\beta_2]$ be a
continuous increasing function and let $\psi: [\beta_1,\beta_2
]\rightarrow [\alpha_1,\alpha_2]$ be its inverse,
$\psi(\phi(\alpha))=\alpha$ (so that $\phi(\alpha_i)=\beta_i$,
$i=1,2$).
Define
\begin{equation}
F(a,b)=\int_{\alpha_1}^{a}\phi \, \dd x+\int_{\beta_1}^{b}\psi \,
\dd x-a b+\alpha_1\beta_1.
\end{equation}
 Young's inequality \cite{hardy34,diaz70,mitrinovic93} states that for every $a \in [\alpha_1, \alpha_2]$ and $b\in [\beta_1,\beta_2]$
\begin{equation}
0\le F(a,b),
\end{equation}
where the equality holds iff $\phi(a)=b$ (or, equivalently,
$\psi(b)=a$).

Among the classical inequalities Young's inequality is probably the
most intuitive. Indeed, its meaning can be  easily grasped once the
integrals are regarded as areas below and on the left of the graph
of $\phi$ (see, for instance, \cite{tolsted64}). Despite its
simplicity, it has profound consequences. For instance, the Cauchy,
Holder and Minkowski inequalities can be easily derived from it
\cite{tolsted64}.

In this work I am going to improve Young's inequality as follows

\begin{theorem} \label{ksd}
Under the  assumptions of Young's inequality, we have  for every $a
\in [\alpha_1,\alpha_2]$ and $b\in [\beta_1,\beta_2]$,
\begin{equation} \label{pdf}
0\le F(a,b)\le -(\psi(b)-a)(\phi(a)-b).
\end{equation}
where the former equality holds if and only if the latter equality
holds.
\end{theorem}

Note that the theorem contains Young's inequality as a special case,
with the advantage that the equality case is naturally taken into
account by the special form of the upper bound. For instance, if
$\psi(b)=a$ then $F(a,b)=0$  which is one of the additional
statements contained in the classical formulation of Young's
inequality. Nevertheless, I will not prove again Young's inequality,
instead I will use  it repeatedly to obtain the  extended version
given by theorem \ref{ksd}.


\begin{remark}
Over the years several extensions of Young's inequality have been
considered. A good account is given by \cite{mitrinovic93}. Among
those only M. Merkle's contribution \cite{merkle74} seems to go in
the same direction considered by this work. Theorem \ref{ksd}
improves Merkle's result, which in the case $\alpha_1=\beta_1=0$
states that (notation of this work)
\[F(a,b) \le \max\{a\phi(a), b\psi(b)\}-ab.\] Indeed,  the last term
of (\ref{pdf}) can be rewritten
\[
[a\phi(a)+b\psi(b)-\phi(a)\psi(b)]-ab,
\]
and we have only to show that
\[
a\phi(a)+b\psi(b)-\phi(a)\psi(b) \le  \max\{a\phi(a), b\psi(b)\},
\]
and that for some $a,b$, the inequality is strict. Indeed, if
$\phi(a)> b$ then, since $\phi$ and $\psi$ are one the inverse of
the other, $a > \psi(b)$ and thus $a\phi(a) > b\psi(b)$. Then
\[
a\phi(a)+b\psi(b)-\phi(a)\psi(b)= a \phi(a)+(b-\phi(a)) \psi(b) <
a\phi(a) = \max\{a\phi(a), b\psi(b)\}
\]
The case $\phi(a)< b$ gives again a strict inequality while the case
$\phi(a)=b$ gives an equality.
\end{remark}


\section{The proof}

The proof of theorem \ref{ksd} is based on the next lemma

\begin{lemma} \label{jhs}
For every $a, \tilde{a} \in [\alpha_1,\alpha_2]$ and $b, \tilde{b}
\in [\beta_1,\beta_2]$, we have
\begin{equation} \label{ouy}
 F(a,b)+F(\tilde{a},\tilde{b})\ge -(\tilde{a}-a)(\tilde{b}-b),
\end{equation}
where the equality holds iff $\tilde{a}=\psi(b)$ and
$\tilde{b}=\phi(a)$.
\end{lemma}

\begin{proof}
Young's inequality gives
\begin{align}
\int_{\alpha_1}^{a}\phi \, \dd x+\int_{\beta_1}^{\tilde{b}} \psi \,
\dd x+\alpha_1\beta_1&\ge a
\tilde{b} \label{j1}\\
\int_{\alpha_1}^{\tilde{a} }\phi \, \dd x+\int_{\beta_1}^{b}\psi \,
\dd x+\alpha_1\beta_1& \ge  \tilde{a} b \label{j2}
\end{align}
then
\begin{align*}
&[\int_{\alpha_1}^{a}\phi \, \dd x+\int_{\beta_1}^{b}\psi \, \dd x-a
b+\alpha_1\beta_1]+[\int_{\alpha_1}^{\tilde{a}}\phi \, \dd
x+\int_{\beta_1}^{\tilde{b}}\psi \, \dd
x-\tilde{a} \tilde{b}+\alpha_1\beta_1]\\
&= [\int_{\alpha_1}^{a}\phi \, \dd x+\int_{\beta_1}^{\tilde{b}}\psi
\, \dd x+\alpha_1\beta_1]+[\int_{\alpha_1}^{\tilde{a}}\phi \, \dd
x+\int_{\beta_1}^{b}\psi \, \dd x+\alpha_1\beta_1]-ab-\tilde{a} \tilde{b} \\
&\ge a \tilde{b} + \tilde{a} b-ab-\tilde{a}
\tilde{b}=-(\tilde{a}-a)(\tilde{b}-b).
\end{align*}
The equality holds iff it holds in  (\ref{j1}) and (\ref{j2}), that
is iff $\tilde{a}=\psi(b)$ and $\tilde{b}=\phi(a)$.
\end{proof}

We are ready to prove the theorem

\begin{proof}[Proof of theorem \ref{ksd}]
Consider  (\ref{ouy})   with $\tilde{a}=\psi(b)$ and
$\tilde{b}=\phi(a)$
\[
F(a,b)+F(\psi(b),\phi(a))= -(\psi(b)-a)(\phi(a)-b).
\]
By Young's inequality, since $\psi(b) \in [\alpha_1,\alpha_2]$ and
$\phi(a) \in [\beta_1,\beta_2]$, $F(\psi(b),\phi(a))\ge 0$, thus
\[
F(a,b)\le -(\phi(a)-b)(\psi(b)-a).
\]
The equality holds iff $F(\psi(b),\phi(a))=0$ which holds, again by
the usual Young's inequality, iff $\phi(\psi(b))=\phi(a)$ i.e.
$b=\phi(a)$ (or equivalently $a=\psi(b)$), which holds iff the
inequality $F(a,b)\ge 0$ is actually an equality.

\end{proof}


\section{The Legendre transform}

It is worthwhile to recall the connection with the Legendre
transform. If $\Phi:[\alpha_1, \alpha_2] \to \mathbb{R} $ and
$\Psi:[\beta_1, \beta_2] \to \mathbb{R} $ are two $C^1$ functions
with increasing derivatives such that  they are  the Legendre
transform of each other then it is well known that they admit the
integral representation $\Phi(a)=\Phi(\alpha_1)+\int_{\alpha_1}^a
\phi \,\dd x$, $\Psi(b)=\Psi(\beta_1)+\int_{\beta_1}^b \psi\, \dd x$
where $\phi$ and $\psi$ are two $C^{0}$ increasing function which
are one the inverse of the other, $\beta_1=\phi(\alpha_1)$ and
$\Phi(\alpha_1)+\Psi(\beta_1)=\alpha_1\beta_1$. Thus the theorem for
the Legendre transforms case takes the following form

\begin{theorem}
If $\Phi:[\alpha_1, \alpha_2] \to \mathbb{R} $ and $\Psi:[\beta_1,
\beta_2] \to \mathbb{R} $ are two $C^1$ functions with increasing
derivatives such that  they are  the Legendre transform of each
other, then for every $a \in [\alpha_1, \alpha_2]$, $b \in [\beta_1,
\beta_2] $
\begin{equation} \label{pdf2}
0\le \Phi(a)+\Psi(b)-a b\le -(\Phi'(a)-b)(\Psi'(b)-a),
\end{equation}
where the former equality holds iff the latter equality holds.
\end{theorem}

\begin{example}
Take $\Phi(a)=\frac{a^{\alpha}}{\alpha}$ and
$\Psi(b)=\frac{b^{\beta}}{\beta}$ with
$\frac{1}{\alpha}+\frac{1}{\beta}=1$, and $\alpha,\beta>1$, then we
obtain the inequalities
\begin{equation}
0\le \frac{a^{\alpha}}{\alpha}+\frac{b^{\beta}}{\beta}-ab\le
-(a^{\alpha-1}-b)(b^{\beta-1}-a),
\end{equation}
in particular the last inequality can be rewritten
\[
b^{\beta-1} a^{\alpha-1}  \le
\frac{1}{\alpha}\,b^{\beta}+\frac{1}{\beta}\,
a^{\alpha}=\frac{1}{\alpha}\,(b^{\beta-1})^{\alpha}+\frac{1}{\beta}\,
(a^{\alpha-1})^{\beta},
\]
\end{example}
that is, it has as expected the same form of Young's inequality.


%
%

%





\section*{Acknowledgments}
This work has been partially supported by GNFM of INDAM and by MIUR
under project PRIN 2005 from Universit\`a di Camerino.

%
%
%

\vspace{1cc}


{\small \noindent Department of Applied Mathematics,\\ Florence
 University,  Via S. Marta 3,  \\
 I-50139 Florence, Italy.\\
E--mail: ettore.minguzzi@unifi.it


}\end{document}